\documentclass[10pt]{article}

\usepackage{amsmath}
\usepackage{amsfonts}
\usepackage{amssymb}
\usepackage[english]{babel}
\usepackage{amsthm}

\def\pf{\par\noindent {\bf Proof}~\par\noindent}
\def\qed{~\hfill{$\square$}\pagebreak[1]\par\medskip\par}

\newcommand{\mR}{\mathbb{R}}
\newcommand{\mC}{\mathbb{C}}
\newcommand{\mN}{\mathbb{N}}

\newcommand{\mZ}{\mathbb{Z}}
\newcommand{\p}{\partial}
\newcommand{\ux}{\underline{x}}

\newcommand{\uom}{\underline{\omega}}
\newcommand{\pux}{\underline{\partial}}

\newcommand{\mcH}{\mathcal{H}}

\newcommand{\mcL}{\mathcal{L}}
\newcommand{\enb}{\overline{e_0}}
\newcommand{\Dbar}{\overline{D}}
\newcommand{\pxn}{\partial_{x_0}}
\newcommand{\onehalf}{\frac{1}{2}}
\newcommand{\la}{\langle}
\newcommand{\ra}{\rangle}

\newtheorem{lemma}{Lemma}
\newtheorem{proposition}{Proposition}
\newtheorem{remark}{Remark}
\newtheorem{corollary}{Corollary}

\numberwithin{equation}{section}
\numberwithin{theorem}{section}
\numberwithin{proposition}{section}
\numberwithin{lemma}{section}
\numberwithin{definition}{section}
\numberwithin{remark}{section}
\numberwithin{corollary}{section}

\begin{document}

\title{Distributional Boundary Values of Harmonic Potentials in Euclidean Half--space as Fundamental Solutions of Convolution Operators in Clifford Analysis}

\author{F.\ Brackx, H.\ De Bie, H.\ De Schepper}

\date{\small{Clifford Research Group, Department of Mathematical Analysis,\\ Faculty of Engineering and Architecture, Ghent University\\
Building S22, Galglaan 2, B-9000 Gent, Belgium\\}}

\maketitle

\begin{abstract}
\noindent In the framework of Clifford analysis, a chain of harmonic and monogenic potentials in the upper half of Euclidean space $\mR^{m+1}$ was recently constructed, including a higher dimensional analogue of the logarithmic function in the complex plane. In this construction the distributional limits of these potentials at the boundary $\mR^{m}$ are crucial. The remarkable relationship between these distributional boundary values and four basic pseudodifferential operators linked with the Dirac and Laplace operators is studied.
\end{abstract}

\maketitle


\section{Introduction}
\label{intro}


In a recent paper \cite{bdbds1} a generalization to Euclidean upper half--space $\mR^{m+1}_+$ was constructed of the  logarithmic function $\ln{z}$ which is holomorphic in the upper half of the complex plane. This construction was carried out in the framework of Clifford analysis, where the functions under consideration take their values in the universal Clifford algebra $\mR_{0,m+1}$ constructed over the Euclidean space $\mR^{m+1}$, equipped with a quadratic form of signature $(0,m+1)$. The concept of a higher dimensional holomorphic function, mostly called monogenic function, is expressed by means of a generalized Cauchy--Riemann operator, which is a combination of the derivative with respect to one of the real variables, say $x_0$, and the so--called Dirac operator $\pux$ in the remaining real variables $(x_1, x_2, \ldots, x_m)$. The generalized Cauchy--Riemann operator $D$ and its Clifford algebra conjugate $\overline{D}$ linearize the Laplace operator, whence Clifford analysis may be seen as a refinement of harmonic analysis.\\[-2mm]

The starting point of the construction of a higher dimensional monogenic logarithmic function, was the fundamental solution of the generalized Cauchy--Riemann operator $D$, also called Cauchy kernel, and its relation to the Poisson kernel and its harmonic conjugate in $\mR^{m+1}_+$. We then proceeded by induction in two directions, {\em downstream} by differentiation and {\em upstream} by primitivation, yielding a doubly infinite chain of monogenic, and thus harmonic, potentials. This chain mimics the well--known sequence of holomorphic potentials in $\mC_+$ (see e.g. \cite{slang}):
$$
\frac{1}{k!} z^k \left[ \ln z - ( 1 + \frac{1}{2} + \ldots + \frac{1}{k}) \right] \rightarrow \ldots \rightarrow z ( \ln z - 1) \rightarrow \ln z 
\stackrel{\frac{d}{dz}}{\longrightarrow} \frac{1}{z} \rightarrow - \frac{1}{z^2} \rightarrow \ldots \rightarrow (-1)^{k-1} \frac{(k-1)!}{z^k}
$$
Identifying the boundary of upper half--space with $\mR^m \cong \{(x_0,\ux) \in \mR^{m+1} : x_0 = 0\}$, the distributional limits for $x_0 \rightarrow 0+$ of those potentials were computed. They split up into two classes of distributions, which are linked by the Hilbert transform, one scalar--valued, the second one Clifford vector--valued. They form two of the four families of Clifford distributions which were thoroughly studied in a series of papers, see \cite{fb1, fb2, distrib} and the references therein.\\

These distributional boundary values are really fundamental, since not only they are used in the definition of the harmonic and monogenic potentials, but also uniquely determine the conjugate harmonic potentials obtained by primitivation, thanks to the simple, but crucial, fact that a monogenic function in $\mR^{m+1}_+$ vanishing at the boundary $\mR^{m}$ indeed is zero. Whence the need to predict the distributional boundary values when constructing the, at that moment unknown, upstream potentials. To that end the distributional boundary values have to be identified  in some way ab initio, which is the aim of the present paper. It is shown that half of them may be recovered as fundamental solutions of specific powers of the Dirac operator, and also half of them, but not the missing ones, as fundamental solutions of specific powers of the Laplace operator. By introducing two new pseudodifferential operators, next to and related to the complex powers of the Dirac and Laplace operators, the whole double infinite set of distributional boundary values may now be identified as fundamental solutions of the four operators. As a remarkable demonstration of symmetry, the distributional boundary values also can serve as convolution kernels for the corresponding pseudodifferential operators of the same kind but with opposite exponent. \\[-2mm]

The organization of the paper is as follows. To make the paper self--contained we recall in Section 2 the basics of Clifford algebra and Clifford analysis and in Section 3 the main results of \cite{bdbds1} on the conjugate harmonic and monogenic potentials in upper half--space $\mR^{m+1}_+$. The four pseudodifferential operators needed for recovering all the distributional boundary values of these harmonic potentials as fundamental solutions, are studied in four consecutive sections. Sections 4 and 6 are devoted to the complex powers of the Dirac and Laplace operator respectively and their fundamental solutions. In Sections 5 and 7 the two new operators, also depending on a complex parameter, and their fundamental solutions are studied. Section 8 contains some conclusions.


\section{Basics of Clifford analysis}
\label{basics}


Clifford analysis (see e.g. \cite{red}) is a function theory which offers a natural and elegant generalization to higher dimension of holomorphic functions in the complex plane and refines harmonic analysis. Let  $(e_0, e_1,\ldots,e_m)$ be the canonical orthonormal basis of Euclidean space $\mR^{m+1}$ equipped with a quadratic form of signature $(0,m+1)$. Then the non--commutative multiplication in the universal real Clifford algebra $\mR_{0,m+1}$ is governed by the rule 
$$
e_{\alpha} e_{\beta} + e_{\beta} e_{\alpha} = -2 \delta_{\alpha \beta}, \qquad \alpha,\beta = 0, 1,\ldots,m
$$
whence $\mR_{0,m+1}$ is generated additively by the elements $e_A = e_{j_1} \ldots e_{j_h}$, where $A=\lbrace j_1,\ldots,j_h \rbrace \subset \lbrace 0,\ldots,m \rbrace$, with $0\leq j_1<j_2<\cdots < j_h \leq m$, and $e_{\emptyset}=1$. 
For an account on Clifford algebra we refer to e.g. \cite{porteous}.\\[-2mm]

We identify the point $(x_0, x_1, \ldots, x_m) \in \mR^{m+1}$ with the Clifford--vector variable 
$$
x = x_0 e_0  + x_1 e_1  + \cdots x_m e_m = x_0 e_0  + \ux
$$ 
and the point $(x_1, \ldots, x_m) \in \mR^{m}$ with the Clifford--vector variable $\ux$. 
The introduction of spherical co--ordinates $\ux = r \uom$, $r = |\ux|$, $\uom \in S^{m-1}$, gives rise to the Clifford--vector valued locally integrable function $\uom$, which is to be seen as the higher dimensional analogue of the {\em signum}--distribution on the real line; we will encounter $\uom$ as one of the distributions discussed below.\\[-2mm]

At the heart of Clifford analysis lies the so--called Dirac operator 
$$
\p = \pxn e_0 + \p_{x_1} e_1 + \cdots \p_{x_m} e_m =  \pxn e_0 + \pux
$$
which squares to the negative Laplace operator: $\p^2 = - \Delta_{m+1}$, while also $\pux^2 = - \Delta_{m}$. The fundamental solution of the Dirac operator $\p$ is given by
$$
E_{m+1} (x) = - \frac{1}{\sigma_{m+1}} \ \frac{x}{|x|^{m+1}}
$$
where $\sigma_{m+1} = \frac{2\pi^{\frac{m+1}{2}}}{\Gamma(\frac{m+1}{2})}$ stands for the area of the unit sphere $S^{m}$ in $\mR^{m+1}$.
We also introduce the generalized Cauchy--Riemann operator 
$$
D = \onehalf \enb \p = \onehalf (\pxn + \enb \pux)
$$ 
which, together with its Clifford algebra conjugate $\Dbar = \onehalf(\pxn - \enb \pux)$, also decomposes the Laplace operator: $D \Dbar = \Dbar D = \frac{1}{4} \Delta_{m+1}$. \\[-2mm]

A continuously differentiable function $F(x)$, defined in an open region $\Omega \subset \mR^{m+1}$ and taking values in the Clifford algebra  $\mR_{0,m+1}$, is called (left--)monogenic if it satisfies  in $\Omega$ the equation $D F = 0$, which is equivalent with $\p F = 0$. \\

Singling out the basis vector $e_0$, we can decompose the real Clifford algebra $\mR_{0,m+1}$ in terms of the Clifford algebra 
$\mR_{0,m}$ as $\mR_{0,m+1} = \mR_{0,m} \oplus \enb \mR_{0,m}$. Similarly we decompose the considered functions as 
$$
F(x_0,\ux) = F_1(x_0,\ux) + \enb F_2(x_0,\ux)
$$ 
where $F_1$ and $F_2$ take their values in the Clifford algebra $\mR_{0,m}$; mimicking functions of a complex variable, we will call $F_1$ the {\em real} part and $F_2$ the {\em imaginary} part of the function $F$.\\[-2mm]

We will extensively use two families of distributions in $\mR^m$, which have been thoroughly studied in \cite{fb1,fb2,distrib}. The first family $\mathcal{T} = \{ T_\lambda : \lambda \in \mC\}$ is very classical. It consists of the radial distributions 
$$
T_\lambda = {\rm Fp} \ r^{\lambda} = {\rm Fp} \ (x_1^2 + \ldots + x_m^2)^{\frac{\lambda}{2}}
$$
their action on a test function $\phi \in \mathcal{S}(\mR^m)$ being given by
$$
\la T_\lambda, \phi \ra = \sigma_m \la {\rm Fp} \; r^\mu_+, \Sigma^{(0)}[\phi] \ra 
$$
with $\mu = \lambda +m-1$. In the above expressions ${\rm Fp}\; r^\mu_+$ stands for the classical "finite part" distribution on the real $r$-axis  and $\Sigma^{(0)}$ is the scalar valued generalized spherical mean, defined on scalar valued test functions $\phi(\ux)$ by
$$
\Sigma^{(0)}[\phi] = \frac{1}{\sigma_m} \int_{S^{m-1}} \phi(\ux) \, dS(\uom)
$$
This family $\mathcal{T}$ contains a.o. the fundamental solutions of the natural powers of the Laplace operator. As convolution operators they give rise to the traditional Riesz potentials (see e.g. \cite{helgason}). The second family $\mathcal{U} = \{ U_\lambda : \lambda \in \mC\}$ of distributions arises in a natural way by the action of the  Dirac operator $\pux$ on $\mathcal{T}$. The $U_{\lambda}$--distributions thus are typical Clifford analysis constructs: they are Clifford--vector valued, and they also arise as products of $T_{\lambda}$--distributions with the distribution $\uom = \frac{\ux}{|\ux|}$, mentioned above. The action of $U_\lambda$ on a test function $\phi \in \mathcal{S}(\mR^m)$ is given by
$$
\la U_\lambda, \phi \ra  = \sigma_m \la {\rm Fp} \; r^\mu_+, \Sigma^{(1)}[\phi] \ra 
$$
with $\mu = \lambda +m-1$, and where the Clifford--vector valued generalized spherical mean $\Sigma^{(1)}$ is defined on scalar valued test functions $\phi(\ux)$ by
$$
\Sigma^{(1)}[\phi] = \frac{1}{\sigma_m} \int_{S^{m-1}} \uom \ \phi(\ux) \, dS(\uom) 
$$
Typical examples in the $\mathcal{U}$--family are the fundamental solutions of the Dirac operator and of its odd natural powers.\\[-2mm]

The normalized distributions $T^{*}_\lambda$ and $U^{*}_\lambda$ arise when removing the singularities of $T_\lambda$ and $U_\lambda$  by dividing them by an appropriate Gamma-function showing the same simple poles. The scalar $T^{*}_\lambda$--distributions are defined by
\begin{equation}
\left \{
\begin{array}{ll}
\displaystyle{T_\lambda^* = \pi^{\frac{\lambda+m}{2}} \frac{T_\lambda}{\Gamma \left ( \frac{\lambda+m}{2} \right )}}, & \lambda \ne -m-2l\\[5mm]
\displaystyle{T_{-m-2l}^* = \frac{\pi^{\frac{m}{2}-l}}{2^{2l} \Gamma \left ( \frac{m}{2} + l \right )} (-\Delta_m)^l \delta (\ux)}, & l \in \mN_0
\label{defTstar}
\end{array}
\right . 
\end{equation}
while the Clifford--vector valued distributions $U^{*}_\lambda$ are defined by
\begin{equation}
\left \{
\begin{array}{ll}
\displaystyle{U_\lambda^* = \pi^{\frac{\lambda+m+1}{2}} \, \frac{U_\lambda}{\Gamma \left ( \frac{\lambda + m + 1}{2} \right )}}, & \lambda \ne -m-2l-1\\[5mm]
\displaystyle{U_{-m-2l-1}^* = - \frac{\pi^{\frac{m}{2}-l}}{2^{2l+1} \, \Gamma \left ( \frac{m}{2} + l + 1 \right )} \; \pux^{2l+1} \delta(\ux)}, & l \in \mN_0
\label{defUstar}
\end{array}
\right . 
\end{equation}
The normalized distributions $T_\lambda^*$ and $U_\lambda^*$ are holomorphic mappings from $\lambda \in \mC$ to the space $\mathcal{S}'(\mR^m)$ of tempered distributions. As already mentioned they are intertwined by the action of the Dirac operator; more generally they enjoy the following properties: for all $\lambda \in \mC$ one has
\begin{itemize}
\item[(i)] $\ux \; T_\lambda^* = \frac{\lambda+m}{2\pi} \; U_{\lambda+1}^*$; \quad
$\ux \; U_\lambda^* = U_\lambda^* \; \ux = - T_{\lambda+1}^*$
\item[(ii)] $\pux \; T_\lambda^* = \lambda \; U_{\lambda-1}^*$; \quad 
$\pux \; U_\lambda^* = U_\lambda^* \; \pux = - 2\pi \; T_{\lambda-1}^*$
\item [(iii)] $\Delta_m T_\lambda^* = 2 \pi \lambda T_{\lambda-2}^*$ ; \quad
$\Delta_m U_\lambda^* = 2 \pi (\lambda-1) U_{\lambda-2}^*$
\item[(iv)] $r^2 T_\lambda^* = \frac{\lambda + m}{2\pi} \ T_{\lambda+2}^*$; \quad
$r^2 U_\lambda^* = \frac{\lambda + m + 1}{2\pi} \ U_{\lambda+2}^*$
\end{itemize}

Of particular importance for the sequel are the convolution formulae for the $T_\lambda^*$-- and $U_\lambda^*$--distributions; we list them in the following proposition and refer the reader to \cite{distrib} for more details. Let us mention that the convolution of the distributions from both families is commutative notwithstanding the Clifford vector character of the  $U_\lambda^*$--distributions.

\begin{proposition}
\label{prop1}
\rule{0mm}{0mm}
\begin{itemize}
\item[(i)] For all $(\alpha,\beta) \in \mC \times \mC$ such that $ \alpha \neq 2j, j \in \mN_0$, $ \beta \neq 2k, k \in \mN_0$ and $\alpha+\beta+m \neq 2l, l \in \mN_0$ the convolution $T_{\alpha}^* \ast T_{\beta}^*$ is the tempered distribution given by 
$$
T_{\alpha}^* \ast T_{\beta}^* = \pi^{\frac{m}{2}}\; \frac{\Gamma \left( -\frac{\alpha + \beta + m}{2} \right)}{\Gamma \left( -\frac{\alpha}{2} \right) \Gamma \left( -\frac{\beta}{2} \right)} \; T_{\alpha+\beta+m}^*
$$

\item[(ii)] For $(\alpha,\beta) \in \mC \times \mC$ such that $ \alpha \neq 2j+1, \  \beta \neq 2k, \ \alpha+\beta \neq -m+2l+1, \ j, k, l \in \mN_0$ one has
$$
U_{\alpha}^*  \ast   T_{\beta}^*  =   T_{\beta}^* \ast  U_{\alpha}^* =  \pi^{\frac{m}{2}}\; \frac{\Gamma \left( -\frac{\alpha + \beta + m - 1}{2} \right)}{\Gamma \left( -\frac{\alpha-1}{2} \right) \Gamma \left( -\frac{\beta}{2} \right)} \; U_{\alpha+\beta+m}^*
$$

\item[(iii)] For $(\alpha,\beta) \in \mC \times \mC$ such that $ \alpha \neq 2j+1, \  \beta \neq 2k+1, \ \alpha+\beta \neq -m+2l, \ j, k, l \in \mN_0$ one has
$$
U_{\alpha}^*  \ast   U_{\beta}^*	  =  U_{\beta}^*  \ast   U_{\alpha}^*  = \pi^{\frac{m}{2}+1} \displaystyle{ \frac{\Gamma(- \frac{\alpha+\beta+m}{2})}{\Gamma(\frac{-\alpha+1}{2}) \Gamma(\frac{-\beta+1}{2})}}  \;  T_{\alpha+\beta+m}^*
$$
\end{itemize}
\end{proposition}


\section{Harmonic and monogenic potentials in $\mR_+^{m+1}$}
\label{potentials}


In this section we gather the most important results on harmonic and monogenic potentials in upper half--space $\mR^{m+1}$, which were established in \cite{bdbds1}.\\[-2mm]

The starting point is the Cauchy kernel of Clifford analysis, i.e. the fundamental solution of the generalized Cauchy--Riemann operator $D$:
$$
C_{-1}(x_0,\ux) = \frac{1}{\sigma_{m+1}} \, \frac{x \overline{e_0}}{|x|^{m+1}} =  \frac{1}{\sigma_{m+1}} \, \frac{x_0 - 
\overline{e_0} \ux}{|x|^{m+1}}
$$
which may be decomposed in terms of the traditional Poisson kernels in $\mR^{m+1}_+$:
$$
C_{-1}(x_0,\ux) = \frac{1}{2} A_{-1}(x_0,\ux) + \frac{1}{2} \overline{e_0} \, B_{-1}(x_0,\ux)
$$
where, also mentioning the usual notations, for $x_0 >0$,
\begin{eqnarray*}
A_{-1}(x_0,\ux) & = & P(x_0,\ux) \ = \ \phantom{-} \frac{2}{\sigma_{m+1}} \, \frac{x_0}{|x|^{m+1}} \label{A-1}\\
B_{-1}(x_0,\ux) & = & Q(x_0,\ux) \ = \ - \frac{2}{\sigma_{m+1}} \, \frac{\ux}{|x|^{m+1}} \label{B-1}
\end{eqnarray*}
Their distributional limits for $x_0 \rightarrow 0+$ are given by
\begin{eqnarray*}
a_{-1}(\ux) & = & \lim_{x_0 \rightarrow 0+} A_{-1}(x_0,\ux) \ = \ \delta(\ux) \ \ = \ \phantom{-} \frac{2}{\sigma_m} \, T^\ast_{-m} \\
b_{-1}(\ux) & = & \lim_{x_0 \rightarrow 0+} B_{-1}(x_0,\ux) \ = \ H(\ux) \ = \ - \frac{2}{\sigma_{m+1}} \, U^\ast_{-m} 
\end{eqnarray*}
where the distribution
$$
H(\ux) = - \frac{2}{\sigma_{m+1}} \, U^\ast_{-m} = - \frac{2}{\sigma_{m+1}} \, \mbox{Pv} \frac{\ux}{|\ux|^{m+1}}
$$
with $\mbox{Pv}$ standing for the "principal value" distribution in $\mR^m$, is the convolution kernel of the Hilbert transform $\mathcal{H}$ in $\mR^m$ (see e.g.\ \cite{gilmur}). Note also that both distributional boundary values are linked by this Hilbert transform:
\begin{eqnarray*}
\mathcal{H} \left [ a_{-1} \right ] & = & \mathcal{H} \left [ \delta \right ] \ = \ H \ast \delta \ = \ H \ = \ b_{-1} \\
\mathcal{H} \left [ b_{-1} \right ] & = & \mathcal{H} \left [ H \right ] \ = \ H \ast H \ = \ \delta \ = \ a_{-1}
\end{eqnarray*}
since $\mathcal{H}^2 = \mathbf{1}$.\\

The first in the sequence of so--called {\em downstream} potentials is the function $C_{-2}$ defined by
$$
\overline{D} C_{-1} = C_{-2} = \frac{1}{2} A_{-2} + \frac{1}{2} \overline{e_0} B_{-2}
$$
Clearly it is monogenic  in $\mR^{m+1}_+$, since $DC_{-2} = D \overline{D} C_{-1} = \frac{1}{4} \Delta_{m+1} C_{-1} = 0$. 
The definition itself of $C_{-2}(x_0,\ux)$ implies that it shows the monogenic potential (or primitive) $C_{-1}(x_0,\ux)$ and the conjugate harmonic potentials $A_{-2}(x_0,\ux)$ and $\overline{e_0} B_{-2}(x_0,\ux)$. The distributional limits for $x_0 \rightarrow 0+$ of these harmonic potentials are given by
$$
\left \{ \begin{array}{rcl}
a_{-2}(\ux) = \lim_{x_0 \rightarrow 0+} A_{-2}(x_0,\ux) & = & \displaystyle\frac{2}{\sigma_{m+1}} \, {\rm Fp} \displaystyle\frac{1}{|\ux|^{m+1}} \ = \ - \displaystyle\frac{4 \pi}{\sigma_{m+1}} T^\ast_{-m-1}\\[4mm]
b_{-2}(\ux) = \lim_{x_0 \rightarrow 0+} B_{-2}(x_0,\ux) & = & - \pux \delta \ = \ \displaystyle\frac{2m}{\sigma_m} \, U^\ast_{-m-1}
\end{array} \right .
$$

Proceeding in the same manner, the sequence of {\em downstream} monogenic potentials in $\mR_+^{m+1}$ is defined by
$$
C_{-k-1} = \overline{D} C_{-k} = \overline{D}^2 C_{-k+1} = \ldots = \overline{D}^k C_{-1}, \qquad k=1,2,\ldots
$$
where each monogenic potential decomposes into two conjugate harmonic potentials:
$$
C_{-k-1} = \frac{1}{2} A_{-k-1} + \frac{1}{2} \overline{e_0} B_{-k-1}, \qquad k=1,2,\ldots
$$
with, for $k$ odd, say $k=2\ell-1$,
$$
\left \{ \begin{array}{rcl}
A_{-2\ell} & = & \p_{x_0}^{2\ell-1} A_{-1} \ = \ - \p_{x_0}^{2\ell-2} \pux B_{-1} \ = \ \ldots \ = \ - \pux^{2\ell-1} B_{-1} \\[2mm]
B_{-2\ell} & = & \p_{x_0}^{2\ell-1} B_{-1} \ = \ - \p_{x_0}^{2\ell-2} \pux A_{-1} \ = \ \ldots \ = \ - \pux^{2\ell-1} A_{-1} 
\end{array} \right .
$$
while for $k$ even, say $k=2\ell$,
$$
\left \{ \begin{array}{rcl}
A_{-2\ell-1} & = & \p_{x_0}^{2\ell} A_{-1} \ = \ - \p_{x_0}^{2\ell-1} \pux B_{-1} \ = \ \ldots \ = \ \pux^{2\ell} A_{-1} \\[2mm]
B_{-2\ell-1} & = & \p_{x_0}^{2\ell} B_{-1} \ = \ - \p_{x_0}^{2\ell-1} \pux A_{-1} \ = \ \ldots \ = \ \pux^{2\ell} B_{-1} 
\end{array} \right .
$$
Their distributional limits for $x_0 \rightarrow 0+$ are given by
$$
\left \{ \begin{array}{rcl}
a_{-2\ell} & = & (- \pux)^{2\ell-1} H  = - 2^{2\ell-1} \displaystyle\frac{\Gamma \left ( \frac{m+2\ell-1}{2} \right )} {\pi^{\frac{m-2\ell+1}{2}}} \, T^{\ast}_{-m-2\ell+1}  \\[5mm]
& = & (-1)^{\ell-1}  2^{\ell-1} (2\ell-1)!!  \displaystyle\frac{\Gamma \left ( \frac{m+2\ell-1}{2} \right )} {\pi^{\frac{m+1}{2}}} \, {\rm Fp} \displaystyle\frac{1}{r^{m+2\ell-1}}                                \\[7mm]
b_{-2\ell} & = & (- \pux)^{2\ell-1} \delta = 2^{2\ell-1} \displaystyle\frac{\Gamma \left (\frac{m+2\ell}{2} \right )} {\pi^{\frac{m-2\ell+2}{2}}} \, U^{\ast}_{-m-2\ell+1}       \end{array} \right .
$$
and
$$
\left \{ \begin{array}{rcl}
a_{-2\ell-1} & = &  \pux^{2\ell} \delta = 2^{2\ell} \displaystyle\frac{\Gamma \left (\frac{m+2\ell}{2} \right )} {\pi^\frac{m-2\ell}{2}}  \, T^{\ast}_{-m-2\ell}                    
\\[7mm]
b_{-2\ell-1} & = & \pux^{2\ell} H  \ = \ - 2^{2\ell} \displaystyle\frac{\Gamma \left ( \frac{m+2\ell+1}{2} \right )} {\pi^\frac{m-2\ell+1}{2}} \, U^{\ast}_{-m-2\ell}  \\[5mm]
& = &(-1)^{\ell-1}  2^{\ell} (2\ell-1)!!  \displaystyle\frac{\Gamma \left ( \frac{m+2\ell+1}{2} \right )}{\pi^\frac{m+1}{2}} \,  {\rm Fp} \displaystyle\frac{1}{r^{m+2\ell}} \, \omega  \end{array} \right .
$$
They show the following properties.
\begin{lemma}
\label{lem2}
One has for $j,k=1,2,\ldots$
\begin{itemize}
\item[(i)] $a_{-k} \xrightarrow{\hspace*{1mm} -\pux \hspace*{1mm}} b_{-k-1} \xrightarrow{\hspace*{1mm} -\pux \hspace*{1mm}} a_{-k-2}$
\item[(ii)] $\mathcal{H} \left [ a_{-k} \right ] = b_{-k}$, $\mathcal{H} \left [ b_{-k} \right ] = a_{-k}$
\item[(iii)] $a_{-j} \ast a_{-k} = a_{-j-k+1}$ \\
$a_{-j} \ast b_{-k} = b_{-j} \ast a_{-k} = b_{-j-k+1}$ \\
$b_{-j} \ast b_{-k} = a_{-j-k+1}$.
\end{itemize}
\end{lemma}

Let us have a look at the so--called {\em upstream} potentials. To start with the fundamental solution of the Laplace operator $\Delta_{m+1}$ in $\mR^{m+1}$, sometimes called Green's function, and here denoted by $\frac{1}{2}A_0(x_0,\ux)$, is given by
$$
\frac{1}{2}A_0(x_0,\ux) = - \frac{1}{m-1} \frac{1}{\sigma_{m+1}} \frac{1}{|x|^{m-1}}
$$
Its conjugate harmonic in $\mR^{m+1}_+$, in the sense of \cite{red}, is
\begin{equation}
B_0(x_0,\ux) = \frac{2}{\sigma_{m+1}} \, \frac{\ux}{|\ux|^m} \, F_m \left ( \frac{|\ux|}{x_0} \right )
\label{B0}
\end{equation}
where
$$
F_m(v) = \int_0^v \frac{\eta^{m-1}}{(1+\eta^2)^\frac{m+1}{2}} \, d\eta = \frac{v^m}{m} \,  _2F_1 \left ( \frac{m}{2},\frac{m+1}{2};\frac{m}{2}+1;-v^2 \right )
$$
with $_2F_1$ a standard hypergeometric function (see e.g. \cite{grad}). Taking into account that
$$
F_m(+\infty) = \int_0^{+\infty} \frac{\eta^{m-1}}{(1+\eta^2)^\frac{m+1}{2}} \, d\eta = \frac{\sqrt{\pi}}{2} \frac{\Gamma \left ( \frac{m}{2} \right )}{\Gamma \left ( \frac{m+1}{2} \right )}
$$
expression (\ref{B0}) leads to the following distributional limit
$$
b_0(\ux) = \lim_{x_0 \rightarrow 0+} B_0(x_0,\ux) = \frac{1}{\sigma_m} \frac{\ux}{|\ux|^m} = \frac{1}{\pi} \, \frac{1}{\sigma_m} \, U^\ast_{-m+1}
$$
while $A_0(x_0,\ux)$ itself shows the distributional limit
$$
a_0(\ux) = \lim_{x_0 \rightarrow 0+} A_0(x_0,\ux) = - \frac{2}{m-1} \, \frac{1}{\sigma_{m+1}} \, {\rm Fp} \frac{1}{|\ux|^{m-1}} = 
- \frac{2}{m-1} \, \frac{1}{\sigma_{m+1}} \, T_{-m+1}^\ast
$$
It is readily seen that $\overline{D} A_0 =  \overline{D} \overline{e_0} B_0 =  C_{-1}$. So $A_0(x_0,\ux)$ and $\overline{e_0} B_0(x_0,\ux)$ are conjugate harmonic potentials (or primitives), with respect to the operator $\overline{D}$, of the Cauchy kernel $C_{-1}(x_0,\ux)$ in $\mR^{m+1}_+$. Putting $C_0(x_0,\ux) = \frac{1}{2} A_0(x_0,\ux) + \frac{1}{2} \overline{e_0} B_0 (x_0,\ux)$, it follows that also $\overline{D} C_0 (x_0,\ux)  = C_{-1}(x_0,\ux)$,
which implies that $C_0(x_0,\ux)$ is a monogenic potential (or primitive) of the Cauchy kernel $C_{-1}(x_0,\ux)$ in $\mR^{m+1}_+$. 
Their distributional boundary values are intimately related, as shown in the following lemma.
\begin{lemma}
\label{lemintiem}
One has
\begin{itemize}
\item[(i)] $-\pux a_0 = b_{-1} = H$; \quad $-\pux b_0 = a_{-1} = \delta$
\item[(ii)] $\mathcal{H} \left [a_0 \right ] = b_0$; \quad $\mathcal{H} \left [b_0 \right ] = a_0$
\end{itemize}
\end{lemma}

\begin{remark}
In the upper half of the complex plane the function $\ln(z)$ is a holomorphic potential (or primitive) of the Cauchy kernel $\frac{1}{z}$ and its real and imaginary components are the fundamental solution $\ln |z|$ of the Laplace operator, and its conjugate harmonic $i\, {\rm arg}(z)$ respectively. By similarity we could say that $C_0(x_0,\ux) = \frac{1}{2} A_0(x_0,\ux) + \frac{1}{2} \overline{e_0} B_0(x_0,\ux)$, being a monogenic potential of the Cauchy kernel $C_{-1}(x_0,\ux)$ and the sum of the fundamental solution $A_0(x_0,\ux)$ of the Laplace operator and its conjugate harmonic $\overline{e_0} B_0(x_0,\ux)$, is a {\em monogenic logarithmic function} in the upper half--space $\mR^{m+1}_+$.
\end{remark}

Inspired by the above mentioned properties, the construction of the sequence of {\em upstream} harmonic and monogenic potentials in $\mR^{m+1}_+$ is continued as follows. Putting
$$
\left \{ \begin{array}{rclcl}
A_1(x_0,\ux) & = & a_0(\cdot) \ast A_0(x_0,\cdot)(\ux) & = & b_0(\cdot) \ast B_0(x_0,\cdot) \\[2mm]
B_1(x_0,\ux) & = & a_0(\cdot) \ast B_0(x_0,\cdot)(\ux) & = & b_0(\cdot) \ast A_0(x_0,\cdot)
\end{array} \right .
$$
it is verified that $\overline{D} A_{-1} =  \overline{D} \overline{e_0} B_{-1} =  C_{0}$, whence $A_1(x_0,\ux)$ and $B_1(x_0,\ux)$ are conjugate harmonic potentials in $\mR^{m+1}_+$ of the function $C_0(x_0,\ux)$. It then follows at once that
$$
C_1(x_0,\ux) = \frac{1}{2} A_1(x_0,\ux) + \frac{1}{2} \overline{e_0} B_1(x_0,\ux)
$$
is a monogenic potential in $\mR^{m+1}_+$ of $C_0$. The distributional limits for $x_0 \rightarrow 0+$ of the conjugate harmonic potentials $A_1$ and $B_1$ are given by
$$
\left \{ \begin{array}{rcl}
a_1(\ux) & = & \lim_{x_0 \rightarrow 0+} A_1(x_0,\ux) \ = \ a_0(\cdot) \ast a_0(\cdot)(\ux) \ = \ b_0(\cdot) \ast b_0(\cdot) (\ux) \\[2mm]
b_1(\ux) & = & \lim_{x_0 \rightarrow 0+} B_1(x_0,\ux) \ = \ a_0(\cdot) \ast b_0(\cdot)(\ux) \ = \ b_0(\cdot) \ast a_0(\cdot) (\ux)
\end{array} \right .
$$
Making use of the calculation rules for the convolution of the $T^\ast$-- and $U^\ast$--distributions (see Section 2, Proposition 2.1), these distributional boundary values are explicitly given by
$$
\left \{ \begin{array}{rclcl}
a_1(\ux) & = & \phantom{-} \displaystyle\frac{1}{\pi} \displaystyle\frac{1}{\sigma_m} \displaystyle\frac{1}{m-2} \, T^\ast_{-m+2} & = & \phantom{-} \displaystyle\frac{1}{\sigma_m} \displaystyle\frac{1}{m-2} \displaystyle\frac{1}{|\ux|^{m-2}}
 \\[5mm]
b_1(\ux) & = & - \displaystyle\frac{1}{\pi} \displaystyle\frac{1}{\sigma_{m+1}} \displaystyle\frac{1}{m-1} \, U^\ast_{-m+2} & = & - \displaystyle\frac{1}{\sigma_{m+1}} \displaystyle\frac{2}{m-1} \displaystyle\frac{\ux}{|\ux|^{m-1}}
\end{array} \right .
$$
They show the following properties.
\begin{lemma}
\label{lem54}
\rule{0mm}{0mm}
\begin{itemize}
\item[(i)] $- \pux a_1 = b_0$, $-\pux b_1 = a_0$
\item[(ii)] $\mathcal{H} \left [ a_1 \right ] = b_1$, $\mathcal{H} \left [ b_1 \right ] = a_1$
\end{itemize}
\end{lemma}

The conjugate harmonic potentials $A_1(x_0,\ux)$ and $B_1(x_0,\ux)$ have been determined explicitly: 
$$
\left \{ \begin{array}{rcl}
A_1(x_0,\ux) &=& \displaystyle\frac{2}{m-1} \, \displaystyle\frac{1}{\sigma_{m+1}} \, \displaystyle\frac{1}{|\ux|^{m-2}} \, F_{m-2} \left ( \displaystyle\frac{\ux}{x_0} \right ) \\[5mm]
B_1(x_0,\ux) &=& \displaystyle\frac{2}{\sigma_{m+1}} \, \displaystyle\frac{x_0 \ux}{|\ux|^m} \, F_m \left ( \displaystyle\frac{|\ux|}{x_0} \right ) - \displaystyle\frac{2}{\sigma_{m+1}} \, \displaystyle\frac{1}{m-1} \, \displaystyle\frac{\ux}{|x|^{m-1}}
\end{array} \right .
$$
Proceeding in a similar way, it is verified that the functions $A_2(x_0,\ux)$ and $B_2(x_0,\ux)$ defined by
$$
\left \{
\begin{array}{rcl}
A_2(x_0,\ux) & = & a_0(\cdot) \ast A_1(x_0,\cdot)(\ux) \ = \ b_0(\cdot) \ast B_1(x_0,\cdot)(\ux) \\[2mm]
B_2(x_0,\ux) & = & a_0(\cdot) \ast B_1(x_0,\cdot)(\ux) \ = \ b_0(\cdot) \ast A_1(x_0,\cdot)(\ux)
\end{array}
\right .
$$
are conjugate harmonic potentials in $\mR^{m+1}_+$ of the function $C_1(x_0,\ux)$. It follows that
$$
C_2(x_0,\ux) = \frac{1}{2} A_2(x_0,\ux) + \frac{1}{2} \overline{e_0} B_2(x_0,\ux)
$$
is a monogenic potential in $\mR^{m+1}_+$ of $C_1$.The distributional limits for $x_0 \rightarrow 0+$ are given by
$$
\left \{ \begin{array}{rcl}
a_2(\ux) & = & \lim_{x_0 \rightarrow 0+} A_2(x_0,\ux) \ = \ a_0 \ast a_1(\ux) \ = \ b_0 \ast b_1(\ux) \\[2mm]
b_2(\ux) & = & \lim_{x_0 \rightarrow 0+} B_2(x_0,\ux) \ = \ a_0 \ast b_1(\ux) \ = \ b_0 \ast a_1(\ux) 
\end{array} \right .
$$
which may be calculated explicitly to be
$$
a_2(\ux) =  - \frac{1}{\pi} \frac{1}{(m-1)(m-3)} \frac{1}{\sigma_{m+1}} \, T^\ast_{-m+3}
$$
and
$$
b_2(\ux) =   \frac{1}{2 \pi^2} \frac{1}{\sigma_{m}} \frac{1}{m-2} U^\ast_{-m+3}
$$
They show the following properties.
\begin{lemma}
\label{lem56}
\rule{0mm}{0mm}
\begin{itemize}
\item[(i)] $-\pux a_2 = b_1$, $- \pux b_2 = a_1$
\item[(ii)] $\mathcal{H} \left [ a_2 \right ] = b_2$, $\mathcal{H} \left [ b_2 \right ] = a_2$
\end{itemize}
\end{lemma}

\noindent The conjugate harmonic potentials $A_2(x_0,\ux)$ and $B_2(x_0,\ux)$ were also explicitly determined:
\begin{eqnarray*}
A_2(x_0,\ux) & = & \frac{2}{m-1} \frac{1}{\sigma_{m+1}} \frac{x_0}{|\ux|^{m-2}} \, F_{m-2} \left ( \frac{|\ux|}{x_0} \right ) - \frac{2}{m-1} \frac{1}{m-3} \frac{1}{\sigma_{m+1}} \frac{1}{|x|^{m-3}} \\
B_2(x_0,\ux) & = & \frac{1}{\sigma_{m+1}} \frac{\ux |x|^2}{|\ux|^{m}} \, F_{m} \left ( \frac{|\ux|}{x_0} \right ) - \frac{m-3}{m-1} \frac{1}{\sigma_{m+1}} \frac{\ux}{|\ux|^{m-2}} \,   F_{m-2} \left ( \frac{|\ux|}{x_0} \right )
\end{eqnarray*}

For general $k=1,2,3,\ldots$ , the following functions in $\mR^{m+1}_+$ are defined recursively, the convolutions being taken in the variable $\ux \in \mR^{m}$:
\begin{eqnarray*}
A_k(x_0,\ux) & = & a_0 \ast A_{k-1} \ = \ a_1 \ast A_{k-2} \ = \ \ldots \ = \ a_{k-1} \ast A_0 \\
& = & b_0 \ast B_{k-1} \ = \ b_1 \ast B_{k-2} \ = \ \ldots \ = \ b_{k-1} \ast B_0 \\
B_k(x_0,\ux) & = & a_0 \ast B_{k-1} \ = \ a_1 \ast B_{k-2} \ = \ \ldots \ = \ a_{k-1} \ast B_0 \\
& = & b_0 \ast A_{k-1} \ = \ b_1 \ast A_{k-2} \ = \ \ldots \ = \ b_{k-1} \ast A_0 \\
\end{eqnarray*}
and
$$
C_k(x_0,\ux) = \frac{1}{2} A_k(x_0,\ux) + \frac{1}{2} \overline{e_0} B_k(x_0,\ux)
$$
It may be verified that $A_k(x_0,\ux)$ and $B_k(x_0,\ux)$ are conjugate harmonic potentials of $C_{k-1}(x_0,\ux)$ in $\mR^{m+1}_+$, while $C_k(x_0,\ux)$ is a monogenic potential of the same $C_{k-1}(x_0,\ux)$ in $\mR^{m+1}_+$. Their distributional boundary values for $x_0 \rightarrow 0+$ are given by the recurrence relations
\begin{eqnarray*}
a_k(\ux) & = & a_0 \ast a_{k-1}\ = \ a_1 \ast a_{k-2} \ = \ \ldots \ = \  a_{k-1} \ast a_0 \\
& = & b_0 \ast b_{k-1}\ = \ b_1 \ast b_{k-2} \ = \ \ldots \ = \  b_{k-1} \ast b_0 \\
b_k(\ux) & = & a_0 \ast b_{k-1}\ = \ a_1 \ast b_{k-2} \ = \ \ldots \ = \  a_{k-1} \ast b_0 \\
& = & b_0 \ast a_{k-1}\ = \ b_1 \ast a_{k-2} \ = \ \ldots \ = \  b_{k-1} \ast a_0 
\end{eqnarray*}
for which the following explicit formulae my be deduced:
$$
\left \{ \begin{array}{rcl}
a_{2k} & =  & -  \displaystyle\frac{1}{2^{2k+1}} \,  \displaystyle\frac{\Gamma(\frac{m-2k-1}{2})}{\pi^{\frac{m+2k+1}{2}}} \; T^{*}_{-m+2k+1} \\[5mm]
a_{2k-1} & =  & \phantom{-} \displaystyle\frac{1}{2^{2k}} \, \displaystyle\frac{\Gamma(\frac{m-2k}{2})}{\pi^{\frac{m+2k}{2}}} \; T^{*}_{-m+2k}
\end{array} \right .
$$
$$
\left \{ \begin{array}{rcl}
b_{2k} & = & \phantom{-}  \displaystyle\frac{1}{2^{2k+1}}  \, \displaystyle\frac{\Gamma(\frac{m-2k}{2})}{\pi^{\frac{m+2k+2}{2}}} \; U^{*}_{-m+2k+1} \\[5mm]
b_{2k-1} & = &  -  \displaystyle\frac{1}{2^{2k}} \,  \displaystyle\frac{\Gamma(\frac{m-2k+1}{2})}{\pi^{\frac{m+2k+1}{2}}} \; U^{*}_{-m+2k} \\[5mm]
\end{array} \right .
$$
These distributional limits show the following properties.
\begin{lemma}
\label{lem58}
One has for $k=1,2,\ldots$:
\begin{itemize}
\item[(i)] $- \pux a_k = b_{k-1}$; $- \pux b_k = a_{k-1}$
\item[(ii)] $\mathcal{H} \left [ a_k \right ] = b_{-1} \ast a_k = b_k$; $\mathcal{H} \left [ b_k \right ] = b_{-1} \ast b_k = a_k$
\end{itemize}
\end{lemma}


\section{Powers of the Dirac operator}
\label{powersdirac}


The complex power of the Dirac operator $\pux$ was already introduced in \cite{dss} and further studied in \cite{distrib}. It is a convolution operator defined by
\begin{eqnarray}
	\pux^{\mu}[ \, . \, ]  = \pux^{\mu}\delta \ast [ \, . \, ] & = & \left[\frac{1+e^{i\pi\mu}}{2} \, \frac{2^{\mu}\Gamma\left(\frac{m+\mu}{2}\right)}{\pi^{\frac{m-\mu}{2}}}\; T^{*}_{-m-\mu} - \frac{1-e^{i\pi\mu}}{2} \, \frac{2^{\mu}\Gamma\left(\frac{m+\mu+1}{2}\right)}{\pi^{\frac{m-\mu+1}{2}}}\; U^{*}_{-m-\mu}\right] \ast [ \, . \, ] \nonumber \\[2mm]
	& = & \frac{2^{\mu}}{\pi^{\frac{m}{2}}} \; \mbox{Fp} \, \frac{1}{\vert\ux\vert^{\mu+m}}\left[\frac{1+e^{i\pi\mu}}{2} \, \frac{\Gamma\left(\frac{m+\mu}{2}\right)}{\Gamma\left(-\frac{\mu}{2}\right)}-\frac{1-e^{i\pi\mu}}{2} \, \frac{\Gamma\left(\frac{m+\mu+1}{2}\right)}{\Gamma\left(-\frac{\mu-1}{2}\right)}\uom\right] \ast [ \, . \, ]
\label{defdiracmu}
\end{eqnarray}
	
In particular for integer values of the  parameter $\mu$, the convolution kernel  $\pux^{\mu}\delta$ is given by
\begin{equation}
\left \{ \begin{array}{rcl}
\pux^{2k}\delta & =  & \phantom{-} \displaystyle\frac{2^{2k}\Gamma\left(\frac{m+2k}{2}\right)}{\pi^{\frac{m-2k}{2}}}\; T^{*}_{-m-2k}\\[5mm]
\pux^{2k+1}\delta & = &  - \displaystyle\frac{2^{2k+1}\Gamma\left(\frac{m+2k+2}{2}\right)}{\pi^{\frac{m-2k}{2}}}\; U^{*}_{-m-2k-1}
\end{array} \right . 
\label{intpowdirac}
\end{equation}
Note that for $k \in \mN_0$ the above expressions (\ref{intpowdirac}) are in accordance with the definitions (\ref{defTstar}) and (\ref{defUstar}). Moreover, if the dimension $m$ is odd, also all negative integer powers of the Dirac operator are defined by (\ref{intpowdirac}). However, if the dimension $m$ is even, the expressions (\ref{intpowdirac}) are no longer valid for $k=-\frac{m}{2}-n$, with $n=0,1,2,\ldots$ in the case of $\pux^{2k}\delta$ and $n=1,2,\ldots$ in the case of $\pux^{2k+1}\delta$. Summarizing, $\pux^\mu$ is defined for all $\mu \in \mC$, except for $\mu = -m, -m-1, -m-2, \ldots$ when $m$ is even. We will define $\pux^\mu$ for those exceptional parameter values further on. First we prove the following fundamental property.

\begin{proposition}
\label{prop41}
For  $\mu, \nu \in \mC$ when $m$ is odd or for $\mu, \nu \in \mC$ such that $\mu$, $\nu$ and $\mu+\nu$ are different from $ - m, - m - 1, - m - 2, \ldots $ when $m$ is even, one has
$$\pux^\mu\delta \ast \pux^\nu\delta =   \pux^{\mu+\nu}\delta$$
\end{proposition}

\pf
Using definition (\ref{defdiracmu}) for $\pux^\mu \delta$ and $\pux^\nu \delta$, the convolution at the left--hand side decomposes into four terms.
They are respectively given by
$$
\frac{1+e^{i\pi\mu}}{2}  \, \frac{1+e^{i\pi\nu}}{2} \, 2^{\mu+\nu} \, \pi^{\frac{m}{2}} \, \frac{  \Gamma \left(\frac{m+\mu+\nu}{2}\right)}{\pi^{\frac{m-\mu}{2}} \, \pi^{\frac{m-\nu}{2}}}\; T^{*}_{-m-\mu-\nu}
$$
for the first one, 
$$
- \frac{1+e^{i\pi\mu}}{2} \, \frac{1-e^{i\pi\nu}}{2} \, 2^{\mu+\nu} \, \pi^{\frac{m}{2}} \, \frac{  \Gamma \left(\frac{m+\mu+\nu+1}{2}\right)}{\pi^{\frac{m-\mu}{2}} \, \pi^{\frac{m-\nu+1}{2}}} \; U^{*}_{-m-\mu-\nu}
$$
for the second,
$$
- \frac{1-e^{i\pi\mu}}{2} \, \frac{1+e^{i\pi\nu}}{2} \, 2^{\mu+\nu} \, \pi^{\frac{m}{2}} \, \frac{  \Gamma\left(\frac{m+\mu+\nu+1}{2}\right)}{\pi^{\frac{m-\mu+1}{2}} \, \pi^{\frac{m-\nu}{2}}}\; U^{*}_{-m-\mu-\nu}
$$
for the third, and
$$
\frac{1-e^{i\pi\mu}}{2} \, \frac{1-e^{i\pi\nu}}{2} \, 2^{\mu+\nu} \, \pi^{\frac{m}{2}} \, \frac{  \Gamma\left(\frac{m+\mu+\nu}{2}\right)}{\pi^{\frac{m-\mu}{2}} \, \pi^{\frac{m-\nu}{2}}} \; T^{*}_{-m-\mu-\nu}
$$
for the fourth. The sum of the first and the fourth term thus equals
$$
\frac{1+e^{i\pi(\mu+\nu)}}{2} \, 2^{\mu+\nu} \, \frac{  \Gamma\left(\frac{m+\mu+\nu}{2}\right)}{\pi^{\frac{m-\mu-\nu}{2}}}\; T^{*}_{-m-\mu-\nu}
$$
while the sum of the second and the third term equals
$$
-\frac{1-e^{i\pi(\mu+\nu)}}{2} \, 2^{\mu+\nu} \, \frac{  \Gamma\left(\frac{m+\mu+\nu+1}{2}\right)}{\pi^{\frac{m-\mu-\nu+1}{2}}}\; U^{*}_{-m-\mu-\nu}
$$
The sum of the latter two expressions is exactly $\pux^{\mu+\nu} \delta$.
\qed

\begin{corollary}
\label{cor41}
For $\mu \in \mC$ when $m$ is odd or for $\mu \in \mC \backslash \{ \pm m, \pm m \pm 1, \pm m \pm 2, \ldots \}$ when $m$ is even, one has
$$\pux^\mu \delta \ast \pux^{-\mu} \delta =   \delta$$
\end{corollary}

Now we put for $\mu \in \mC$ when $m$ is odd or for $\mu \in \mC \backslash \{ m, m + 1, m + 2, \ldots \}$ when $m$ is even
$$
E_\mu = \pux^{-\mu} \delta=  \frac{1+e^{-i\pi\mu}}{2} \frac{2^{-\mu}\Gamma\left(\frac{m-\mu}{2}\right)}{\pi^{\frac{m+\mu}{2}}}\; T^{*}_{-m+\mu} - \frac{1-e^{-i\pi\mu}}{2} \frac{2^{-\mu}\Gamma\left(\frac{m-\mu+1}{2}\right)}{\pi^{\frac{m+\mu+1}{2}}}\; U^{*}_{-m+\mu}
$$
and in particular for $k \in \mZ$ when $m$ is odd or for $k \in \mZ \backslash \{\frac{m}{2}+n, n=0,1,2,\ldots  \}$ when $m$ is even
$$
\left \{ \begin{array}{rcl}
E_{2k} & =  & \phantom{-}  \displaystyle\frac{1}{2^{2k}} \displaystyle\frac{\Gamma\left(\frac{m-2k}{2}\right)}{\pi^{\frac{m+2k}{2}}}\; T^{*}_{-m+2k}\\[5mm]
E_{2k+1} & = &  - \displaystyle\frac{1}{2^{2k+1}} \displaystyle\frac{\Gamma\left(\frac{m-2k}{2}\right)}{\pi^{\frac{m+2k+2}{2}}}\; U^{*}_{-m+2k+1}
\end{array} \right .
$$
Then Corollary \ref{cor41} implies that, for $\mu \in \mC$ when $m$ is odd or for $\mu \in \mC \backslash \{ \pm m, \pm m \pm 1, \pm m \pm 2, \ldots \}$ when $m$ is even, $E_\mu = \pux^{-\mu} \delta$ is the fundamental solution of the operator $\pux^\mu$:
$$
\pux^\mu E_\mu = \pux^\mu \delta \ast E_\mu = \delta
$$
This is in accordance with a result in \cite{distrib}.\\[-2mm]

It is also clear that, in the case where the dimension $m$ is even, once the fundamental solutions $E_{m+n}$ of $\pux^{m+n}, n=0,1,2,\ldots$ are known, we can use these expressions for defining the operators $\pux^{-m-n}, n=0,1,2,\ldots$. 
To that end we recall a result of \cite{distrib}.

\begin{proposition}
\label{prop42}
If the dimension $m$ is even, for $n=0,1,2,\ldots$, the fundamental solution $E_{m+n}$ of the operator $\pux^{m+n}$ is given by
$$
\left \{ \begin{array}{rcl}
E_{m+2j} & =  &  (p_{2j} \ln{r} + q_{2j}) \; T^{*}_{2j}\\[1mm]
E_{m+2j+1} & = & (p_{2j+1} \ln{r} + q_{2j+1}) \; U^{*}_{2j+1}
\end{array} \right . \qquad j=0,1,2,\ldots
$$
where the constants $p_n$ and $q_n$ satisfy the recurrence relations
$$
\left \{ \begin{array}{rcl}
p_{2j+2} & =  &  \displaystyle\frac{1}{2j+2} \, p_{2j+1}\\[4mm]
q_{2j+2} & = & \displaystyle\frac{1}{2j+2} \, (q_{2j+1} - \displaystyle\frac{1}{2j+2} \, p_{2j+1})
\end{array} \right . \qquad j=0,1,2,\ldots 
$$
and
$$ 
\left \{ \begin{array}{rcl}
p_{2j+1} & =  & - \displaystyle\frac{1}{2\pi} \, p_{2j}\\[4mm]
q_{2j+1} & = & - \displaystyle\frac{1}{2\pi} \, (q_{2j} - \displaystyle\frac{1}{m+2j} \, p_{2j})
\end{array} \right . \qquad j=0,1,2,\ldots 
$$
with starting values $p_{0} =  - \displaystyle\frac{1}{2^{m-1}\pi^m}$ and $q_{0} = 0$.
\end{proposition}

\noindent Now putting, for $m$ even and $n=0,1,2,\ldots$, $\pux^{-m-n} \delta = E_{m+n}$, and hence
$$
\pux^{-m-n} [ \, . \, ] = \pux^{-m-n} \delta \ast [ \, . \, ] = E_{m+n} \ast [ \, . \, ]
$$
we indeed have
$$
\pux^{-m-n} E_{-m-n} = \pux^{-m-n}  \delta \ast \pux^{m+n} \delta = E_{m+n} \ast \pux^{m+n} \delta = \delta
$$
So the operator $\pux^{\mu} [ \, . \, ]$ eventually is defined for all $\mu \in \mC$, and there holds in distributional sense
$$
\pux^\mu [E_\mu]  =  \pux^\mu [\pux^{-\mu} \delta] = \delta , \quad \mu \in \mC
$$
or, at the level of the operators: $\pux^\mu \pux^{-\mu} = {\bf 1}$.


\section{A new operator}
\label{new1}


Recalling the following distributional boundary values  of the conjugate harmonic potentials studied in \cite{bdbds1}
$$
\left \{ \begin{array}{rcl}
a_{2k-1} & =  & \displaystyle\frac{1}{2^{2k}} \, \displaystyle\frac{\Gamma(\frac{m-2k}{2})}{\pi^{\frac{m+2k}{2}}} \; T^{*}_{-m+2k}, \quad k \in \mZ, \ \ \ 2k < m\\[4mm]
b_{2k} & = &   \displaystyle\frac{1}{2^{2k+1}}  \, \displaystyle\frac{\Gamma(\frac{m-2k}{2})}{\pi^{\frac{m+2k+2}{2}}} \; U^{*}_{-m+2k+1}, \quad k \in \mZ, \ \ \ 2k < m
\end{array} \right .
$$
it becomes clear, in view of the results in Section \ref{powersdirac}, that these distributional boundary values are nothing but fundamental solutions of appropriate integer powers of the Dirac operator. We have indeed, for integer $k$ such that $2k < m$, that
$$
\left \{ \begin{array}{rclcl}
a_{2k-1} & =  & \phantom{-} E_{2k} & = & \phantom{-} \pux^{-2k} \delta\\[2mm]
b_{2k} & = & - E_{2k+1} & = & - \pux^{-2k-1} \delta
\end{array} \right .
$$
showing that at the same time they are also distributions resulting from the action of the opposite integer powers of the Dirac operator on the delta distribution.\\[-2mm]

It also becomes clear that the other distributional boundary values $a_{2k}$ and $b_{2k-1}$ cannot be expressed in a similar way as fundamental solutions of integer powers of the Dirac operator. Whence the need for a new operator, depending upon a complex parameter $\mu$, the fundamental solutions of which correspond to those distributional boundary values $a_{2k}$ and $b_{2k-1}$. Taking into account the Hilbert pair relationship between the distributional boundary values, we define to that end the operator $^\mu \mcH$ by
$$
^\mu \mcH [ \, . \, ] = \pux^\mu H \ast [ \, . \, ]
$$
where the convolution kernel $\pux^\mu H$ is given by
$$
\pux^{\mu} H  =  \frac{1-e^{i\pi\mu}}{2} \, \frac{2^{\mu}\Gamma\left(\frac{m+\mu}{2}\right)}{\pi^{\frac{m-\mu}{2}}} \; T^{*}_{-m-\mu} - \frac{1+e^{i\pi\mu}}{2} \, \frac{2^{\mu}\Gamma\left(\frac{m+\mu+1}{2}\right)}{\pi^{\frac{m-\mu+1}{2}}} \; U^{*}_{-m-\mu}
$$
The notation for this new kernel is motivated by the fact that, as shown by a straightforward calculation, it may indeed be obtained as $\pux^{\mu} H  = \pux^{\mu} \delta \ast H$. In particular for integer values of the parameter $\mu$, the convolution kernel $\pux^{\mu} H $ reduces to
\begin{equation}
\left \{ \begin{array}{rcl}
\pux^{2k}H & =  & - 2^{2k}  \displaystyle\frac{\Gamma(\frac{m+2k+1}{2})}{\pi^{\frac{m-2k+1}{2}}} \; U^{*}_{-m-2k} \\[5mm]
\pux^{2k+1}H & = &  2^{2k+1}  \displaystyle\frac{\Gamma(\frac{m+2k+1}{2})}{\pi^{\frac{m-2k-1}{2}}} \; T^{*}_{-m-2k-1} 
\end{array} \right .
\label{intpowhilbertdirac}
\end{equation}
with $2k \neq -m-1, -m-3, \ldots$ when $m$ is odd. Note that for $\mu=0$ the operator $^0 \mcH$ reduces to the Hilbert transform, while for $\mu=1$ the so--called Hilbert--Dirac operator (see \cite{fbhds}) is obtained:
$$
^1 \mcH[ \, . \, ] = (-\Delta_m)^\onehalf [ \, . \, ] = \pux H \ast [ \, . \, ] = 2 \frac{\Gamma(\frac{m+1}{2})}{\pi^{\frac{m-1}{2}}} \; T^{*}_{-m-1} \ast [ \, . \, ]
$$
More generally, we also have for integer $k$ such that $2k \neq -m-1, -m-3, \ldots$ when $m$ is odd, 
$$
^{2k+1} \mcH [ \, . \, ] = \pux^{2k+1} H \ast [ \, . \, ] =  (-\Delta_m)^{k + \onehalf} [ \, . \, ] 
$$
Summarizing, the operator $^\mu \mcH$ is defined for all complex values of the parameter $\mu$ except for $\mu = -m, -m-1, -m-2,\ldots$ when $m$ is odd. We will use the same method as above, via the fundamental solutions, to define $^\mu \mcH$ for those exceptional values. 

\begin{proposition}
\label{prop51}
For  $\mu, \nu \in \mC$ when $m$ is even or for $\mu, \nu \in \mC$ such that $\mu$, $\nu$ and $\mu+\nu$ are different from $ - m, - m - 1, - m - 2, \ldots $ when $m$ is odd, one has
$$
\pux^\mu H \ast \pux^\nu H =   \pux^{\mu+\nu} H
$$
\end{proposition}

\pf
The proof is similar to that of Proposition \ref{prop41}. \qed

\begin{corollary}
\label{cor51}
For $\mu \in \mC$ when $m$ is even or for $\mu \in \mC \backslash \{ \pm m, \pm m \pm 1, \pm m \pm 2, \ldots \}$ when $m$ is odd, one has
$$\pux^\mu H \ast \pux^{-\mu} H =   \delta$$
\end{corollary}

Now we put for $\mu \in \mC$ when $m$ is even or for $\mu \in \mC \backslash \{ m, m + 1, m + 2, \ldots \}$ when $m$ is odd
$$
F_\mu = \pux^{-\mu} H =  \frac{1- e^{-i\pi\mu}}{2} \frac{2^{-\mu}\Gamma\left(\frac{m-\mu}{2}\right)}{\pi^{\frac{m+\mu}{2}}}\; T^{*}_{-m+\mu} - \frac{1+ e^{-i\pi\mu}}{2} \frac{2^{-\mu}\Gamma\left(\frac{m-\mu+1}{2}\right)}{\pi^{\frac{m+\mu+1}{2}}}\; U^{*}_{-m+\mu}
$$
and in particular for integer values of the parameter $\mu$
$$
\left \{ \begin{array}{rcl}
F_{2k} & =  & - \displaystyle\frac{1}{2^{2k}} \displaystyle\frac{\Gamma\left(\frac{m-2k+1}{2}\right)}{\pi^{\frac{m+2k+1}{2}}}\; U^{*}_{-m+2k}, \ \ k \neq \frac{m+1}{2}+n, n=0,1,\ldots {\rm when} \ m \ {\rm is \ odd}\\[5mm]
F_{2k+1} & = &   \phantom{-} \displaystyle\frac{1}{2^{2k+1}} \displaystyle\frac{\Gamma\left(\frac{m-2k-1}{2}\right)}{\pi^{\frac{m+2k+1}{2}}}\; T^{*}_{-m+2k+1}, \ \ k \neq \frac{m-1}{2}+n, n=0,1,\ldots {\rm when} \ m \ {\rm is \ odd}
\end{array} \right .
$$

Then Corollary \ref{cor51} implies that for $\mu \in \mC$ when $m$ is even or for $\mu \in \mC \backslash \{ \pm m, \pm m \pm 1, \pm m \pm 2, \ldots\}$ when $m$ is odd
$$
^\mu \mcH [F_\mu] = \pux^\mu H \ast F_\mu = \delta 
$$
expressing the fact that $F_\mu = \pux^{-\mu} H$ is the fundamental solution of the operator $\pux^\mu \mcH$ for the allowed values of $\mu$. So it becomes clear that, in the case where the dimension $m$ is odd, if we succeed in establishing the fundamental solutions $F_{m+n}, n=0,1,2,\ldots$ of the corresponding operators $^{m+n}\mcH$, we can use these expressions for defining the operators $^{-m-n}\mcH$. We first prove that $E_\mu$ and $F_\mu$ form a Hilbert pair.

\begin{proposition}
\label{prop52}
For $\mu \in \mC \backslash \{ m, m + 1, m + 2, \ldots \}$ one has
$$ \mcH [E_\mu]  =   F_\mu$$
\end{proposition}

\pf
For the allowed values of $\mu$ we consecutively have
$$
\mcH  [ E_\mu ] =  ^0\mcH [E_\mu]  =  H \ast E_\mu =  H \ast \pux^{-\mu} \delta =  \pux^{-\mu} \delta \ast H = \pux^{-\mu} H = F_\mu
$$
\qed

Now we determine the fundamental solutions $F_{m+n}, n=0,1,2,\ldots$ when the dimension $m$ is odd.
The general expression for $\pux^{-\mu}$ being no longer valid in that case, this needs a specific approach, which is similar to the one used for determining the fundamental solutions $E_{m+n}$ of $\pux^{-m-n}$ when $m$ was even.

\begin{proposition}
If the dimension $m$ is odd, then, for $n=0,1,2,\ldots$, the fundamental solution of $^{m+n}\mcH$ is given by
$$
\left \{ \begin{array}{rcl}
F_{m+2j} & =  &  (p_{2j} \ln{r} + q_{2j})  \; T^{*}_{2j} \\[2mm]
F_{m+2j+1} & = & ( p_{2j+1} \ln{r} + q_{2j+1})    \; U^{*}_{2j+1} 
\end{array} \right . \qquad j=0,1,2,\ldots
$$
with the same constants $(p_n,q_n)$ as in Proposition \ref{prop42}
\end{proposition}

\pf
We have to prove that $^{m+n}\mcH [F_{m+n}] = \delta$ or $\pux^{m+n}H \ast F_{m+n} = \delta$, or still $\pux^{m+n} \ast F_{m+n} = H$, which will be satisfied if $\pux F_{m+n} = F_{m+n-1}$. If $n$ is even, say $n=2j$, we have
$$
\pux F_{m+2j} = p_{2j} \frac{\ux}{r^2} \; T^*_{2j} + (p_{2j} \ln{r} + q_{2j})  \; \pux T^{*}_{2j} = p_{2j} U^*_{2j-1} + 2j  (p_{2j} \ln{r} + q_{2j}) \; U^*_{2j-1}
$$
from which it follows that the following recurrence relations should hold
$$
\left \{ \begin{array}{rcl}
p_{2j} + 2j q_{2j} & =  &  q_{2j-1}\\[2mm]
2j p_{2j} & = & p_{2j-1}
\end{array} \right .
\quad {\rm or} \quad
\left \{ \begin{array}{rcl}
p_{2j} & =  &  \displaystyle\frac{1}{2j} \, p_{2j-1}\\[4mm]
q_{2j} & = &  \displaystyle\frac{1}{2j} \, (q_{2j-1} - \frac{1}{2j} \, p_{2j-1})
\end{array} \right .
$$
If $n$ is odd, say $n=2j+1$, we have
\begin{eqnarray*}
\pux F_{m+2j+1} &=& p_{2j+1} \frac{\ux}{r^2} \; U^*_{2j+1} + (p_{2j+1} \ln{r} + q_{2j+1})  \; \pux U^{*}_{2j+1} \\
&=&  - p_{2j+1} \frac{2\pi}{m+2j} \; T^*_{2j} - 2\pi  (p_{2j+1} \ln{r} + q_{2j+1}) \; T^*_{2j}
\end{eqnarray*}
leading to the recurrence relations
$$
\left \{ \begin{array}{rcl}
-  \displaystyle\frac{2\pi}{m+2j} p_{2j+1} - 2\pi q_{2j+1} & =  &  q_{2j}\\[3mm]
- 2 \pi p_{2j+1} & = & p_{2j}
\end{array} \right .
\quad {\rm or} \quad
\left \{ \begin{array}{rcl}
p_{2j+1} & =  & -  \displaystyle\frac{1}{2\pi} \, p_{2j}\\[3mm]
q_{2j+1} & = &  -  \displaystyle\frac{1}{2\pi} \, (q_{2j} -  \displaystyle\frac{1}{m+2j} \, p_{2j})
\end{array} \right .
$$
\qed

\noindent So putting for $m$ odd and $n=0,1,2,\ldots$, $\pux^{-m-n} H = F_{m+n}$, and hence
$$
^{-m-n} \mcH [ \, . \, ] = \pux^{-m-n} H \ast [ \, . \, ] = F_{m+n} \ast [ \, . \, ]
$$
we indeed have
$$
^{-m-n} \mcH [F_{-m-n}]  = \pux^{-m-n} H \ast \pux^{m+n} H  = \delta
$$
Eventually the operator $^{\mu} \mcH$ is defined for all $\mu \in \mC$, and there holds in distributional sense
$$
^{\mu} \mcH  [F_{\mu}]  \ = \ ^{\mu} \mcH  [^{-\mu} H] = \delta
$$
or, at the level of operators: $^{\mu} \mcH  ^{-\mu} \mcH  = {\bf 1}$.\\[-2mm]

Now we expect the distributional boundary values 
$$
\left \{ \begin{array}{rclll}
a_{2k} & =  & -  \displaystyle\frac{1}{2^{2k+1}} \,  \displaystyle\frac{\Gamma(\frac{m-2k-1}{2})}{\pi^{\frac{m+2k+1}{2}}} \; T^{*}_{-m+2k+1}, & k \in \mZ, & 2k+1 < m\\[3mm]
b_{2k-1} & = &  -  \displaystyle\frac{1}{2^{2k}} \,  \displaystyle\frac{\Gamma(\frac{m-2k+1}{2})}{\pi^{\frac{m+2k+1}{2}}} \; U^{*}_{-m+2k}, & k \in \mZ, & 2k-1 < m
\end{array} \right .
$$
to be fundamental solutions of $^{\mu} \mcH$ for specific values of $\mu$. This is indeed the case since
$$
\left \{ \begin{array}{rclclll}
a_{2k} & =  & - F_{2k+1} & = & - \pux^{-2k-1} H, & k \in \mZ, & 2k+1 < m\\[2mm]
b_{2k-1} & = &  \phantom{-} F_{2k} &  = & \phantom{-} \pux^{-2k} H, & k \in \mZ, & 2k-1 < m
\end{array} \right .
$$
We conclude that all distributional boundary values of the sequence of conjugate harmonic potentials of Section \ref{potentials} are fundamental solutions of $\pux^\mu$ and $^\mu \mcH$ for specific integer values of $\mu$.


\section{Powers of the Laplace operator}
\label{powerslaplace}


For complex powers of the Laplace operator the standard definition (see \cite{helgason}) reads $(-\Delta_m)^\beta [ \, . \, ] = (-\Delta_m)^\beta \delta \ast [ \, . \, ]$, where the convolution kernel $(-\Delta_m)^\beta \delta$ is given by
\begin{eqnarray*}
\label{deflaplacebeta}
	(-\Delta_m)^{\beta}\delta & = &  2^{2\beta} \, \frac{ \Gamma\left(\frac{m+2\beta}{2}\right)}{\pi^{\frac{m-2\beta}{2}}}\; T^{*}_{-m-2\beta}
\end{eqnarray*}
Whence apparently $(-\Delta_m)^\beta$ is defined for all complex values of the parameter $\beta$, except for $\beta = -\frac{m}{2}, -\frac{m}{2}-1, -\frac{m}{2}-2, \ldots$.\\[-2mm]
	
In particular for integer values $k$ of the parameter $\beta$, except for $k= -\frac{m}{2}, -\frac{m+2}{2}, -\frac{m+4}{2}, \ldots$ in case the dimension $m$ is even, we have, also in view of \ref{intpowdirac},
\begin{equation}
(-\Delta_m)^{k}\delta  =   2^{2k} \, \frac{ \Gamma\left(\frac{m+2k}{2}\right)}{\pi^{\frac{m-2k}{2}}}\; T^{*}_{-m-2k} = \pux^{2k}
\delta
\label{intlapl}
\end{equation}
which is in accordance with the factorization of the Laplace operator by the Dirac operator. In case the dimension $m$ is odd, we have, also in view of \ref{intpowhilbertdirac},
\begin{equation}
(-\Delta_m)^{k+\onehalf}\delta  =   2^{2k+1} \, \frac{ \Gamma\left(\frac{m+2k+1}{2}\right)}{\pi^{\frac{m-2k-1}{2}}}\; T^{*}_{-m-2k-1} = \pux^{2k+1} H
\label{intonehalflapl}
\end{equation}
for integer $k$, except for $ -\frac{m+1}{2}, -\frac{m+3}{2}, -\frac{m+5}{2}, \ldots$.\\[-2mm]

By a straightforward calculation, similar to the one in the proof of Proposition \ref{prop41}, the following fundamental property is proven.

\begin{proposition}
\label{prop61}
For $\alpha, \beta \in \mC$ such that $\alpha$, $\beta$ and $\alpha+\beta$ are different from $- \frac{m}{2} - n, n=0,1,2,\ldots$ one has
$$
(-\Delta_m)^\alpha \delta \ast   (-\Delta_m)^\beta \delta =   (-\Delta_m)^{\alpha + \beta} \delta
$$
\end{proposition}

\begin{corollary}
\label{cor61}
For $\beta \in \mC \backslash \{ \pm \frac{m}{2}, \pm \frac{m+2}{2}, \pm \frac{m+4}{2}, \ldots \}$ one has
$$
(-\Delta_m)^\beta \delta \ast   (-\Delta_m)^{-\beta} \delta =   \delta
$$
\end{corollary}

Now putting for $\beta \in \mC \backslash \{ \frac{m}{2}, \frac{m+2}{2}, \frac{m+4}{2}, \ldots \}$
$$
K_\beta = (-\Delta_m)^{-\beta} \delta = 2^{-\beta} \, \frac{ \Gamma\left(\frac{m-2\beta}{2}\right)}{\pi^{\frac{m+2\beta}{2}}}\; T^{*}_{-m+2\beta}
$$
and in particular for integer $k$
\begin{equation}
\left \{ \begin{array}{rclclll}
K_{k} & = & (-\Delta_m)^{-k} \delta & = &  \displaystyle\frac{1}{2^{2k}} \, \displaystyle\frac{\Gamma\left(\frac{m-2k}{2}\right)}{\pi^{\frac{m+2k}{2}}} \; T^{*}_{-m+2k}, & 2k < m \ {\rm when} \ m \ {\rm is \ even}\\[4mm]
K_{k+\onehalf} & = & (-\Delta_m)^{-k-\onehalf} \delta & = &   \displaystyle\frac{1}{2^{2k+1}} \, \displaystyle\frac{\Gamma\left(\frac{m-2k-1}{2}\right)}{\pi^{\frac{m+2k+1}{2}}} \; T^{*}_{-m+2k+1}, & 2k < m-1 \ {\rm when} \ m \ {\rm is \ odd}
\end{array} \right .
\label{K}
\end{equation}
the above Corollary \ref{cor61} implies that
$$
(-\Delta_m)^\beta [K_\beta] = \delta, \quad \beta \in \mC \backslash \{ \pm \frac{m}{2}, \pm \frac{m+2}{2}, \pm \frac{m+4}{2}, \ldots \}
$$
which expresses the fact that $K_\beta = (-\Delta_m)^{-\beta} \delta$ is the fundamental solution of the operator $(-\Delta_m)^\beta$ for $\beta \in \mC \backslash \{ \pm \frac{m}{2}, \pm \frac{m+2}{2}, \pm \frac{m+4}{2}, \ldots \}$.\\[-2mm]

We still need to define the operator $(-\Delta_m)^\beta$ for $\beta = - \frac{m}{2}, - \frac{m+2}{2}, - \frac{m+4}{2}, \ldots $ and the fundamental solution $K_\beta$ for $\beta =  \frac{m}{2},  \frac{m+2}{2},  \frac{m+4}{2}, \ldots$. Keeping in mind the formulae
(\ref{intlapl}) and (\ref{intonehalflapl}), which we still want to remain valid, we put, for $n=0,1,2,\ldots$\\
(i) when $m$ is odd:
$$
(-\Delta_m)^{-\frac{m}{2}-n} \delta = \pux^{-m-2n} H = F_{m+2n} = (p_{2n} \ln{r} + q_{2n}) \, T^*_{2n}
$$
and $K_{\frac{m}{2}+n} = F_{m+2n}$; \\[-2mm]

\noindent (ii) when $m$ is even:
$$
(-\Delta_m)^{-\frac{m}{2}-n} \delta = \pux^{-m-2n} \delta = E_{m+2n} = (p_{2n} \ln{r} + q_{2n}) \, T^*_{2n}
$$
and $K_{\frac{m}{2}+n} = E_{m+2n}$.\\[-2mm]

\noindent We then indeed have\\[1mm]
(i) for $m$ odd
$$
(-\Delta_m)^{\frac{m}{2}+n} \left [ K_{\frac{m}{2}+n} \right ] =  \pux^{m+2n} H \ast K_{\frac{m}{2}+n} =  \pux^{m+2n} H \ast F_{m+2n} = \delta
$$
(ii) for $m$ even
$$
(-\Delta_m)^{\frac{m}{2}+n} \left [ K_{\frac{m}{2}+n} \right ] = \pux^{m+2n} \delta \ast K_{\frac{m}{2}+n} =  \pux^{m+2n} \delta \ast E_{m+2n} = \delta
$$
which eventually leads to
$$
(-\Delta_m)^{\frac{m}{2}+n}   (-\Delta_m)^{-\frac{m}{2}-n}= {\bf 1}
$$

Note that for {\em natural} powers of the Laplace operator, the above fundamental solutions are in accordance with the results of \cite{arons}, where also the closed form of the coefficients $p_{2n}$ and $q_{2n}, n=0,1,2,\ldots$ can be found.


\section{A second new operator}
\label{new2}


The conclusion of Sections \ref{powersdirac} and \ref{new1} was that all distributional boundary values of the sequence of conjugate harmonic potentials studied in \cite{bdbds1}, and recalled in Section \ref{potentials}, are fundamental solutions of the operators $\pux^\mu$ and $^\mu \mcH$ for specific integer values of the parameter $\mu$. Wondering if they are also fundamental solutions of the operator $(-\Delta_m)^\beta$ for some specific values of the complex parameter $\beta$, we indeed find that, in view of (\ref{K}),
$$
\left \{ \begin{array}{rclclcll}
a_{2k} & =  & - \displaystyle\frac{1}{2^{2k+1}} \, \displaystyle\frac{\Gamma(\frac{m-2k-1}{2})}{\pi^{\frac{m+2k+1}{2}}} \; T^{*}_{-m+2k+1} 
& = & - K_{k+\onehalf} & = & - (-\Delta_m)^{-k-\onehalf} \delta, & 2k+1 < m\\[4mm]
a_{2k-1} & = &  \phantom{-} \displaystyle\frac{1}{2^{2k}} \, \displaystyle\frac{\Gamma(\frac{m-2k}{2})}{\pi^{\frac{m+2k}{2}}} \; U^{*}_{-m+2k}
& = &   \phantom{-}  K_{k} & = &  \phantom{-} (-\Delta_m)^{-k} \delta, & 2k < m
\end{array} \right .
$$

To recover the distributional boundary values $b_{2k}$ and $b_{2k-1}$ as fundamental solutions of powers of the Laplace operator, apparently a new operator has to come into play again. Bearing in mind that   the distributional boundary values are forming  Hilbert pairs , we define the operator $^\beta \! \mcL$ by
$$
^\beta \! \mcL [ \, . \, ] = (-\Delta_m)^\beta H \ast [ \, . \, ]
$$
where the convolution kernel $(-\Delta_m)^\beta H$ is given by
$$
(-\Delta_m)^\beta H  = -  2^{\beta} \, \frac{ \Gamma\left(\frac{m+2\beta +1}{2}\right)}{\pi^{\frac{m-2\beta+1}{2}}}\; U^{*}_{-m-2\beta} 
$$
The notation for this second new kernel is motivated by the fact that, as shown by a straightforward calculation, it may indeed be obtained as the convolution $(-\Delta_m)^\beta H  = (-\Delta_m)^\beta \delta  \ast H$. Apparently the operator $^\beta \! \mcL$ is defined for all complex values of the parameter $\beta$ except for $\beta = - \frac{m+1}{2} -n, n=0,1,2,\ldots$.\\[-2mm]

Note the particular cases for integer $k$:
$$
\left \{ \begin{array}{rcll}
(-\Delta_m)^{k}H & =  & - 2^{2k} \, \displaystyle\frac{\Gamma(\frac{m+2k+1}{2})}{\pi^{\frac{m-2k+1}{2}}} \; U^{*}_{-m-2k}, & k \neq - \frac{m+1}{2}, - \frac{m+3}{2}, \ldots (m  \ {\rm odd})\\[4mm]
(-\Delta_m)^{k+\onehalf}H & = &  - 2^{2k+1} \, \displaystyle\frac{\Gamma(\frac{m+2k+2}{2})}{\pi^{\frac{m-2k}{2}}} \; U^{*}_{-m-2k-1}, & \quad k \neq - \frac{m+2}{2}, - \frac{m+4}{2}, \ldots (m \ {\rm even})
\end{array} \right .
$$

\noindent It follows that 
\begin{equation}
\left \{ \begin{array}{rcllll}
(-\Delta_m)^{k}H & =  & \pux^{2k} H, & k \in \mZ, & k \neq - \frac{m+1}{2}, - \frac{m+3}{2}, \ldots & (m  \ {\rm odd})\\[2mm]
(-\Delta_m)^{k+\onehalf}H & = & \pux^{2k+1} \delta, & k \in \mZ, & k \neq - \frac{m+2}{2}, - \frac{m+4}{2}, \ldots & (m \ {\rm even})
\end{array} \right .
\label{intpowhilbertlapl}
\end{equation}
or, at the level of the operators: $^k \!\mcL \ = \ ^{2k}\mcH$ and $^{k+\onehalf} \! \mcL \ = \ \pux^{2k+1}$.

\begin{proposition}
\label{prop71}
For $\alpha, \beta \in \mC$ such that $\alpha$, $\beta$ and $\alpha+\beta$ are different from $- \frac{m+1}{2}-n$, $n=0,1,2,\ldots$, one has
$$
(-\Delta_m)^\alpha H \ast (-\Delta_m)^\beta H = (-\Delta_m)^{\alpha+\beta} \delta
$$
\end{proposition}

\pf
The proof is similar to the one of Proposition \ref{prop41}. \qed

\begin{corollary}
\label{cor71}
For $\beta \in \mC \backslash \{ \pm  \frac{m+1}{2} \pm n, n=0,1,2, \ldots \}$ one has $(-\Delta_m)^\beta H \ast (-\Delta_m)^{-\beta} H =  \delta$.
\end{corollary}

Putting,  for $\beta \in \mC \backslash \{  \frac{m+1}{2} + n, n=0,1,2, \ldots \}$,
$$
L_\beta = (-\Delta_m)^{-\beta} H =  -  2^{-\beta} \, \frac{ \Gamma\left(\frac{m-2\beta +1}{2}\right)}{\pi^{\frac{m+2\beta+1}{2}}}\; U^{*}_{-m+2\beta} 
$$
and in particular for integer $k$
$$
\left \{ \begin{array}{rcll}
L_{k} & =  &  -  2^{-2k} \, \displaystyle\frac{ \Gamma\left(\frac{m-2k +1}{2}\right)}{\pi^{\frac{m+2k+1}{2}}}\; U^{*}_{-m+2k}, & 2k < m+1 \\[3mm]
L_{k+\onehalf} & =  &  -  2^{-2k+1} \, \displaystyle\frac{ \Gamma\left(\frac{m-2k}{2}\right)}{\pi^{\frac{m+2k+2}{2}}}\; U^{*}_{-m+2k+1}, & 2k < m 
\end{array} \right .
$$
Corollary \ref{cor71} implies that
$$
^\beta \! \mcL [L_\beta] = \delta, \quad \beta \in \mC \backslash \{ \pm  \frac{m+1}{2} \pm n, n=0,1,2, \ldots \}
$$
expressing the fact that $L_\beta = (-\Delta_m)^{-\beta} H$ is the fundamental solution of the operator $^\beta \! \mcL$ for $\beta \in \mC \backslash \{ \pm  \frac{m+1}{2} \pm n, n=0,1,2, , \ldots \}$. 

\begin{proposition}
\label{prop72}
For $\beta \in \mC \backslash \{ \frac{m+n}{2}, n=0,1,2,\ldots \}$ one has $ \mcH [K_\beta]  =   L_\beta$.
\end{proposition}

\pf
For the allowed values of $\beta$ we consecutively have
$$
\mcH  [ K_\beta ] =  ^0\mcH [K_\beta]  =  H \ast K_\beta =  H \ast (-\Delta_m)^{-\beta} \delta =  (-\Delta_m)^{-\beta}  H   = L_\beta
$$
\qed

Now we define the operator $^\beta \! \mcL$ for $\beta = - \frac{m+1}{2}-n, n=0,1,2,\ldots$ and the fundamental solution $L_\beta$ for $\beta =  \frac{m+1}{2}+n, n=0,1,2,\ldots$:\\[3mm]
(i) if $m$ is odd, we put
$$
(-\Delta_m)^{ - \frac{m+1}{2}-n} H = \pux^{-m-2n-1} H = F_{m+2n+1}  = (p_{2n+1} \ln{r} + q_{2n+1}) \; U^*_{2n+1}
$$
and $L_{ \frac{m+1}{2}+n} = F_{m+2n+1}$;\\[-2mm]

\noindent (ii) if $m$ is even, we put
$$
(-\Delta_m)^{ - \frac{m+1}{2}-n} H = \pux^{-m-2n-1} \delta = E_{m+2n+1}  = (p_{2n+1} \ln{r} + q_{2n+1}) \;U^*_{2n+1}
$$
and $L_{ \frac{m+1}{2}+n} = E_{m+2n+1}$.\\[-2mm]

\noindent In this way the properties (\ref{intpowhilbertlapl}) are preserved for the exceptional values of $\beta$, and moreover
$$
^{ \frac{m+1}{2}+n} \! \mcL [L_{ \frac{m+1}{2}+n} ] \ = \ ^{\frac{m+1}{2}+n} \! \mcL [(-\Delta_m)^{ - \frac{m+1}{2}-n} H] = \delta
$$
or
$$
\left( ^{ \frac{m+1}{2}+n} \mcL \right) \left( ^{- \frac{m+1}{2}-n} \mcL \right) = {\bf 1}
$$
As expected the distributional boundary values $b_\ell$ are indeed recovered from the fundamental solutions of the operator $^\beta \! \mcL$, since
$$
\left \{ \begin{array}{rclcl}
b_{2k} & =  & - E_{2k+1} & = & - L_{k+\onehalf}, \quad 2k < m\\
b_{2k-1} & =  & \phantom{-} F_{2k} & = & \phantom{-} L_{k}, \quad 2k < m+1
\end{array} \right .
$$


\section{Conclusion}


In this paper we have shown that the distributional boundary values for $x_0  \rightarrow 0+$ of the sequence of conjugate harmonic potentials in upper half--space $\mR_+^{m+1} = \{ x_0 e_0 + \ux : \ux \in \mR^m, x_0 > 0\}$, can be expressed as fundamental solutions of specific powers of four operators: the standard operators $\pux^\mu$ and $(-\Delta_m)^\beta$ and the two newly introduced operators $^\mu \mcH$ and $^\beta \mcL$. The extension of the definition of those four operators to the exceptional values of the complex parameters $\mu$ and $\beta$ for which the operators were not defined initially, was crucial for this purpose. The unifying character of the families of Clifford distributions $\mathcal{T}^\ast$ and $\mathcal{U}^\ast$ in this is remarkable.\\[-2mm]

For specific values of the complex parameters $\mu$ and $\beta$, the four operators studied are interconnected. Since these relationships are fundamental we recall them here. For integer $k$ one has for the corresponding convolution kernels
$$
\begin{array}{rclcrcl}
(-\Delta_m)^k \delta & = & \pux^{2k} \delta & \hspace*{5mm} & (-\Delta_m)^k H & = & \pux^{2k} H\\[5mm]
(-\Delta_m)^{k+\onehalf} \delta & = & \pux^{2k+1} H & \hspace*{5mm} & (-\Delta_m)^{k+\onehalf} H & = & \pux^{2k+1} \delta
\end{array}
$$
The apparent symmetries in these formulae strengthen the idea that, like the Dirac or delta--distribution $\delta$, also the Hilbert kernel $H$ really is a fundamental distribution, more or less a counterpart to the pointly supported $\delta$.\\[-2mm]

The above formulae also generalize the well--known fact that the composition of the two Clifford vector operators $\pux$ and $\mcH$ equals the scalar operator {\em square root of the Laplacian} $(-\Delta_m)^\onehalf$:
$$
(-\Delta_m)^\onehalf \delta = \pux H
$$
This also leads to the well--known scalar factorization of the Laplace operator in terms of pseudodifferential operators: $-\Delta_m = (-\Delta_m)^\onehalf (-\Delta_m)^\onehalf$, next to its vector factorization, used by P.A.M.\ Dirac under matrix disguise: $- \Delta_m = \pux \; \pux$. The fact that the Laplace operator may be factorized in two completely different ways is explained by the fact that both the convolution of two $T^*$--distributions --- to which family $(-\Delta_m)^\onehalf$ belongs --- and the convolution of two $U^*$--distributions --- to which family $\pux$ belongs --- result into a $T^*$--distribution. On the contrary, vector valued operators, such as $\pux$ and $\mcH$, only have one factorization based on the convolution of a $T^*$--distribution with a $U^*$--distribution, as shown in the following remarkable formulae:
$$
\pux \delta = (-\Delta_m)^\onehalf H \quad \quad \quad H = (-\Delta_m)^\onehalf \delta
$$
Finally note that each of the fundamental solutions of the four operators studied in this paper, or, in other words, each of the distributional boundary values of the conjugate harmonic potentials studied in \cite{bdbds1}, may be used as a convolution kernel to define an operator of the same kind but with opposite parameter value:
$$
\begin{array}{rclclcl}
\pux^\mu E_\mu & = & \delta \quad & {\rm and} & \quad E_\mu \ast [ \, . \, ] & = & \pux^{-\mu} [\, . \,]\\[2mm]
\pux^\mu H F_\mu & = & \delta \quad & {\rm and} & \quad F_\mu \ast [\, . \, ] & = & ^{-\mu} \mcH [\, . \,]\\[2mm]
(-\Delta_m)^\beta K_\beta & = & \delta \quad & {\rm and} & \quad K_\beta \ast [\, . \,] & = &  (-\Delta_m)^{-\beta} [\, . \,]\\[2mm]
(-\Delta_m)^\beta H L_\beta & = & \delta \quad & {\rm and} & \quad L_\beta \ast [\, . \,] & = &  ^{-\beta} \mcL [\, . \,]

\end{array}
$$


\end{document}